\documentclass[12pt,a4paper,openright,reqno]{extarticle}
\usepackage[english]{babel}
\usepackage{newlfont}
\usepackage{amsbsy}
\usepackage[dvips]{graphicx}
\usepackage{color}
\usepackage{bbm}
\usepackage{graphicx, subfigure}
\usepackage{color}
\usepackage[center]{caption2}
\usepackage{amssymb}
\usepackage{amsmath}

\usepackage{latexsym}
\usepackage{amsthm}
\usepackage{amsthm}

\usepackage[dvips]{geometry}
\theoremstyle{plain}
\newtheorem{thm}{Theorem}

\newtheorem{lem}[thm]{Lemma}

\newtheorem{nota}[thm]{Notation}

\theoremstyle{remark}
\newtheorem{rem}[thm]{Remark}

\usepackage{mathtools}
\def\alphaltiset#1#2{\ensuremath{\left(\kern-.2em\left(\genfrac{}{}{0pt}{}{#1}{#2}\right)\kern-.2em\right)}}
\usepackage[center]{caption2}

\begin{document}

\title{Affine subspaces of matrices with constant rank }
\author{Elena Rubei}
\date{}
\maketitle

{\footnotesize\em Dipartimento di Matematica e Informatica ``U. Dini'', 
viale Morgagni 67/A,
50134  Firenze, Italia }

{\footnotesize\em
E-mail address: elena.rubei@unifi.it}

\def\thefootnote{}
\footnotetext{ \hspace*{-0.36cm}


{\bf \copyright } This manuscript version is made available under the CC-BY-NC-ND 4.0 license https://creativecommons.org/licenses/by-nc-nd/4.0/}

\begin{abstract} 
For every $m,n \in \mathbb{N}$ and every field $K$, let $M(m \times n, K)$ be the vector space of the $(m \times n)$-matrices over $K$ and let $S(n,K)$ be the vector space of the symmetric $(n \times n)$-matrices over  $K$.
We say that an affine subspace $S$ of $M(m \times n, K)$ or of  $S(n,K)$ has constant rank $r$ if every matrix of $S$ has rank $r$.
Define
$${\cal A}^K(m \times n; r)=  \{ S \;| \; S \; \mbox{\rm  affine subsapce of $M(m \times n, K)$  of constant rank }  r\}$$ $${\cal A}_{sym}^K(n;r)=  \{ S \;| \; S \; \mbox{\rm  affine subsapce of $S(n,K)$  of constant rank }  r\}$$
$$a^K(m \times n;r) = \max \{\dim S \mid S \in {\cal A}^K(m \times n; r ) \}.$$
$$a_{sym}^K(n;r) = \max \{\dim S \mid S \in {\cal A}_{sym}^K(n,r)  \}.$$

In this paper we prove the following two formulas for $r \leq m \leq n$: $$a_{sym}^{\mathbb{R}}(n;r) \leq
\left\lfloor \frac{r}{2} \right\rfloor
\left(n-  \left\lfloor \frac{r}{2} \right\rfloor\right)$$ 
$$a^{\mathbb{R}}(m \times n;r) = r(n-r)+ \frac{r(r-1)}{2} .$$ 

\end{abstract}

\section{Introduction}

For every $m,n \in \mathbb{N}$ and every field $K$, let $M(m \times n, K)$ be the vector space of the $(m \times n)$-matrices over $K$ and let $S(n,K)$ be the vector space of the symmetric $(n \times n)$-matrices over  $K$. Moreover, denote the $\mathbb{R}$-vector space of the hermitian $(n\times n)$-matrices by $H(n)$.

We say that an affine subspace $S$ of $M(m \times n, K)$ or of  $S(n,K)$ (or of $H(n)$) has constant rank $r$ if every matrix of $S$ has rank $r$ and we say that a linear subspace $S$ 
 of $M(m \times n, K)$ or of $S(n,K)$ has constant rank $r$ if every nonzero matrix of $S$ has rank $r$.

Define
$${\cal A}^K(m \times n; r)=  \{ S \;| \; S \; \mbox{\rm  affine subsapce of $M(m \times n, K)$  of constant rank }  r\}$$
$${\cal A}_{sym}^K(n;r)=  \{ S \;| \; S \; \mbox{\rm  affine subsapce of $S(n,K)$  of constant rank }  r\}$$
$${\cal A}_{herm}(n;r)=  \{ S \;| \; S \; \mbox{\rm  affine subsapce of $H(n)$  of constant rank }  r\}$$
$${\cal A}_{sym}^{\mathbb{R}}(n;p, \nu)=  \{ S \;| \; S \; \mbox{\rm  affine subsapce of $S(n,\mathbb{R})$ s.t.
each $ A \in S$  has signature }  (p, \nu)\}$$
$${\cal A}_{herm}(n;p, \nu)=  \{ S \;| \; S \; \mbox{\rm  affine subsapce of $H(n)$ s.t.
each $ A \in S$  has signature }  (p, \nu)\}$$
$${\cal L}^K (m \times n;r) =  \{ S \;| \; S \; \mbox{\rm  linear subsapce of $M(m \times n, K)$   of constant rank }  r\}$$
$${\cal L}_{sym}^K (n; r) =  \{ S \;| \; S \; \mbox{\rm  linear subsapce of $S(n,K)$  of constant rank }  r\}$$

Let 
$$a^K(m \times n;r) = \max \{\dim S \mid S \in {\cal A}^K(m \times n; r ) \}$$
$$a_{sym}^K(n;r) = \max \{\dim S \mid S \in {\cal A}_{sym}^K(n;r)  \}$$
$$a_{herm}(n;r) = \max \{\dim S \mid S \in {\cal A}_{herm}(n;r)  \}$$
$$a_{sym}^{\mathbb{R}}(n;p, \nu) = \max \{\dim S \mid S \in {\cal A}_{sym}^{\mathbb{R}}(n; p,\nu)  \}$$
$$a_{herm}(n;p, \nu) = \max \{\dim S \mid S \in {\cal A}_{herm}(n; p,\nu)  \}$$
$$l^K(m \times n;r) = \max \{\dim S \mid S \in {\cal L}^K (m \times n;r)\}$$
$$l_{sym}^K(n;r) = \max \{\dim S \mid S \in {\cal L}_{sym}^K (n,r)\}.$$

There is a wide literature on linear subspaces of constant rank. In particular we quote the following theorems:

\begin{thm}  {\bf (Westwick, \cite{W1})} For $2 \leq r \leq m \leq n$, we have:
$$ n-r+1 \leq  l^{\mathbb{C}}(m \times n;r) \leq  m+ n -2 r+1$$
\end{thm}

\begin{thm}  {\bf (Ilic-Landsberg, \cite{I-L})} If $r$ is even and greater than or equal to $2$, then
$$l_{sym}^{\mathbb{C}}(n;r) =n-r +1$$
\end{thm}

In case $r$ odd, the following result holds, see \cite{I-L}, \cite{G}, \cite{H-P}:

\begin{thm}  If $r$ is odd, then
$$l_{sym}^{\mathbb{C}}(n;r) =1$$
\end{thm}

We mention also that,  in \cite{Fl}, Flanders
proved that, if $ r \leq m \leq n$, a linear subspace of $M(m \times n ,\mathbb{C})$ such that every of its elements has rank less than or equal to $r$ has dimension less than or equal to $r n$.

In this paper we investigate on the maximal dimension of affine subspaces of constant rank.  
The main  theorems we prove are the following.

\begin{thm}\label{mio}  Let $n,r \in \mathbb{N}$ with $r \leq n$.
Then
$$a_{sym}^{\mathbb{R}}(n;r) \leq
\left\lfloor \frac{r}{2} \right\rfloor
\left(n-  \left\lfloor \frac{r}{2} \right\rfloor\right)
.$$ 
\end{thm}

\begin{thm}\label{mio2}  Let $m,n,r \in \mathbb{N}$ with $r \leq m \leq n$.
Then
$$a^{\mathbb{R}}(m \times n;r) = rn- \frac{r(r+1)}{2} .$$ 
\end{thm}

We prove also a statement 
on  the maximal dimension of affine subspaces with constant signature in the space of symmetric real matrices, see Theorem \ref{signature}, and one 
on the maximal dimension of affine  subspaces of constant rank in the space of the hermitian matrices, see Theorem \ref{herm}.

\section{Proofs of the theorems}

\begin{nota} Let $m,n \in \mathbb{N} -\{0\} $ and $K$ be a field. 

We denote the $n \times n$ identity matrix over $K$ by
$I_n^K$ (or by $I_n$ when the field is clear from the context).

We denote $E_{i,j}^{K,n}$ the $n \times n$  matrix over $K$ such that 
$$ (E_{i,j}^{K,n})_{x,y} = \left\{ \begin{array}{ll}
1 & \mbox{\rm if} \; (x,y)=(i,j) \\
0 & \mbox{\rm otherwise}
\end{array} \right.$$
We omit the superscript when it is clear from the context.

For any $A \in M(m \times n, K) $ we denote the submatrix of $A$ given by the rows $i_1, , \ldots , i_k$ and the columns $j_1, \ldots, j_s$  by $A^{(j_1, \ldots, j_s)}_{(i_1, \ldots, i_k)}$.
\end{nota}

\begin{lem} \label{lemma2}
Let $n \in \mathbb{N}-\{0\}$ and let $A \in S(n,\mathbb{R})$. Then there exists $s \in \mathbb{R}$ such that $\det(I_n+s A) = 0$ if and only if $A \neq 0$.
\end{lem}

\begin{proof}
$\Rightarrow$ This implication is obvious.

$\Leftarrow$ Suppose $A \neq 0$. Then $A$ has a nonzero eigenvalue $\lambda$. Let $s = - \frac{1}{\lambda}$. Then $$\det (I_n+s A)= 
s^n \det \left( \frac{1}{s} I_n +A\right) = s^n \det( A - \lambda I_n) =0.$$

\end{proof}

\begin{lem} \label{lemma3}
Let $r \in \mathbb{N}-\{0\}$. Let $K$ be a field such that, if $x \in K^r -\{0\}$, then $x_1^2+\ldots + x_r^2 \neq 0$.

 Then,  for any   $A \in M(r \times r, K)$ and $x 
\in K^r -\{0\}$, we have that 
$$ \det \begin{pmatrix} I_r+sA & sx\\ s \,{}^t  x & 0 \end{pmatrix} $$
is a nonzero polynomial in $s$.
\end{lem}

\begin{proof} The statement follows immediately from the fact that the coefficient of $s^2$ in 
$ \det \begin{pmatrix} I_r+sA & sx\\ s \,{}^t  x & 0 \end{pmatrix} $ is
 $-(x_1^2 + \ldots+ x_r^2)$.
\end{proof}

\begin{lem} \label{lemma4}
Let $r \in \mathbb{N}-\{0\}$. 
Let $A \in H(r)$ be  positive-definite or negative-definite  and $x \in \mathbb{C}^r-\{0\}$. Then the matrix
$$\begin{pmatrix} A & x\\ {}^t \overline{x} & 0 \end{pmatrix} $$
is invertible.
\end{lem}

\begin{proof} If $A$ is positive-definite,  up to elementary row operations and the same elementary column  operations on the first $r$ rows and the first $r$ columns, we can suppose that $A=I_r$.
If $\begin{pmatrix} x \\ 0\end{pmatrix}$ were linear combination of the first $r$ columns of 
$\begin{pmatrix} I & x \\ {}^t \overline{x} & 0 \end{pmatrix} $, we would have that $0=|x_1|^2 + \ldots+ |x_r|^2$, which is absurd. Analogously if $A$ is negative-definite.
\end{proof}

\begin{rem} \label{abn}
Let $a,b, n \in \mathbb{N}$ with $a+b \leq n$. If $b \geq a $, then $(n-b)b\geq (n-a)a$.
\end{rem}

\begin{proof} Observe that 
$(n-b)b\geq (n-a)a $ if and only if $b^2 - a^2 \leq 
n(b-a) $, which is equivalent to  $b+a \leq n$ (since $b-a \geq 0$), which is true by assumption.
\end{proof}

\begin{proof}[Proof of Theorem \ref{mio}]
 Let $R \in  {\cal A}_{sym}^{\mathbb{R}} (n;r) $. 
 We want to prove that $\dim (R) \leq
 \left\lfloor \frac{r}{2} \right\rfloor
\left(n-  \left\lfloor \frac{r}{2} \right\rfloor\right)
 $. We can write $R$ as $M+ L$ where $M \in S(n, \mathbb{R})$ and $L$ 
  is a linear subspace of $S(n, \mathbb{R})$. Let $Q$ be an invertible matrix such that ${}^t Q M Q$ is a  diagonal matrix $D$ whose diagonal is $(1, \ldots, 1 , -1, \ldots, -1, 0, \ldots, 0)$, where $1$ is repeated $p$ times for some $p$ and $-1$ is repeated $q$ times with $p+q=r$. 
 Let $V= {}^t Q L Q$ and $S=  {}^t Q R Q =D+V$. Obviously $S  \in  {\cal A}_{sym}^{\mathbb{R}} (n;r) $; moreover, $ \dim(S)= \dim (R)$, so to prove that $\dim (R) \leq 
\left\lfloor \frac{r}{2} \right\rfloor
\left(n-  \left\lfloor \frac{r}{2} \right\rfloor\right) 
 $ it is sufficient to prove that $\dim (S) \leq 
 \left\lfloor \frac{r}{2} \right\rfloor
\left(n-  \left\lfloor \frac{r}{2} \right\rfloor\right)$.

 Let $Z$ be the vector subspace of $M(n \times n, \mathbb{R})$ generated by the matrices $E_{i,j}+E_{j,i}$ for $i, j \in \{1, \ldots, p\}$ with $i \leq j$.
 
 Let $U$ be the vector subspace of $M(n \times n, \mathbb{R})$ generated by the matrices $E_{i,j}+E_{j,i}$ for $i, j \in \{p+1, \ldots, r\}$ with $i \leq j$.

 Let $W$ be the vector subspace of $M(n \times n, \mathbb{R})$ generated by the matrices $E_{i,j}+E_{j,i}$ for $i, j \in \{r+1, \ldots, n\}$ with $i \leq j$.

Let $G$    be the vector subspace of $M(n \times n, \mathbb{R})$ generated by the matrices $E_{i,j}+E_{j,i}$ for $i \in \{r+1, \ldots, n\}$, $j\in \{p+1, \ldots , r\}$.

 We want to prove that 
$$V \cap (Z+U+W+G)=\{0\}$$

Let $A \in Z, B \in U , C \in W$, $H \in G$ such that $A+B+C+H \in V$.

\begin{itemize}

\item If there existed $h \in \{r+1, \ldots, n\}$ such that 
$C_{h,h}=0$ and $H_{(h)} \neq 0$, take $s \in \mathbb{R} -\{0\}$ such that  $\det(I_p+sA^{(1,\ldots,p)}_{(1,\ldots,p)}) \neq 0$ and $-I_q+sB^{(p+1, \ldots, r)}_{(p+1, \ldots, r)} $ is negative-definite; then, by Lemma \ref{lemma4}, the matrix
$ \left(
\begin{array}{cc} -I_q+sB^{(p+1, \ldots, r)}_{(p+1, \ldots, r)}  & s \, {}^t \!   (H^{(p+1, \ldots, r)}_{(h)}) \\ 
s H^{(p+1, \ldots, r)}_{(h)} & 0
\end{array} 
\right) $ would be  invertible,  so 
$D+ s(A+B+C+H)$ would have rank greater than $r$,
  so $S $ 
 would not be of constant rank $r$, which is contrary to our assumption.

\item  Suppose there exists $h \in \{r+1, \ldots, n\}$ such that 
$C_{h,h} \neq 0$ and $H_{(h)} \neq 0$; then 
$ \det \left(
\begin{array}{cc}
-I_q+sB^{(p+1, \ldots, r)}_{(p+1, \ldots, r)}  & s \, {}^t \!    (H^{(p+1, \ldots, r)}_{(h)}) \\ 
s H^{(p+1, \ldots, r)}_{(h)} & s \,C_{h,h}
\end{array} 
\right) $ is a polynomial in $s$ with term of degree $1$ equal to $\pm C_{h,h}$, so a nonconstant polynomial. Hence, for $s $ different from a finite number of real numbers, such a determinant is nonzero and then we can find $s$ such that 
$ \det \left(
\begin{array}{cc}
-I_q+sB^{(p+1, \ldots, r)}_{(p+1, \ldots, r)}  & s \, {}^t \!   (H^{(p+1, \ldots, r)}_{(h)}) \\ 
s H^{(p+1, \ldots, r)}_{(h)} & s\,C_{h,h}
\end{array} 
\right) \neq 0 $ and $\det (I_p+sA^{(1, \ldots,p)}_{(1, \ldots,p)}) \neq 0$. So ${\operatorname{rk}} (D+s (A+B+C+H))$ would be greater than $r$, 
  so $S $ 
 would not be of constant rank $r$, which is contrary to our assumption. 
 
 Hence we can conclude that $H=0$.

\item
If  $C$ were nonzero, take $s\in \mathbb{R}$ such that $ \det (I_p +s A^{(1,\ldots,p)}_{(1,\ldots,p)}) \neq 0$ and $ \det (-I_q+ s B^{(p+1, \ldots,r)}_{(p+1, \ldots,r)}) \neq 0$; then $D+s (A+B+C+H)$, that is $D+s (A+B+C)$,  would have rank greater than $r$,  so $S $ 
 would not be of constant rank $r$, which is contrary to our assumption. 
So $C$ must be zero. 

\item
If at least one of $A$ and $B$ were nonzero, take $s \in \mathbb{R}$ such that  
$$\det (I_p+ s A^{(1,\ldots,p)}_{(1,\ldots,p)})=0$$ or 
$$ \det (-I_q+ s B^{(p+1, \ldots,r)}_{(p+1, \ldots,r)}) = 0$$
(there exists by Lemma \ref{lemma2});
 then  $D+ s(A+B+C+H) $, that is 
$D + s (A+B)$, has rank less than $r$, so $S $ 
 would not be of constant rank $r$, which is contrary to our assumption. So also $A$ and $B$ must be zero.
\end{itemize} 

So we have proved that 
$V \cap (Z+U+W+G)=\{0\}$. Hence 
$$ \dim (S) = \dim ( V) \leq  \dim (S(n , \mathbb{R})) - \dim (Z+U+W+G)= p(n-p)$$

In an analogous way we can prove that 
$$ \dim (S) = \dim ( V) \leq  q(n-q).$$

So
\begin{equation} \label{dimSmin}
\dim (S) \leq 
\min \{ p(n-p), q(n-q)\}. 
\end{equation}
 
Observe that 
\begin{equation} \label{minr2}
\min \{ p(n-p), q(n-q)\} \leq
 \left\lfloor \frac{r}{2} \right\rfloor
\left(n-  \left\lfloor \frac{r}{2} \right\rfloor\right),
\end{equation}
in fact: suppose for instance that $p \leq q $, then,
by Remark \ref{abn}, we have that 
\begin{equation} \label{min} 
\min \{ p(n-p), q(n-q)\}= p(n-p);
\end{equation}
moreover, observe that $p \leq q $ ad $p+q=r$ imply that $ p \leq \left\lfloor \frac{r}{2} \right\rfloor$; by applying again Remark \ref{abn} with $(a,b)= \left(p, \left\lfloor \frac{r}{2} \right\rfloor\right)$, we get
\begin{equation} \label{pr2} 
p(n-p) \leq 
 \left\lfloor \frac{r}{2} \right\rfloor
\left(n-  \left\lfloor \frac{r}{2} \right\rfloor\right);
\end{equation}
from (\ref{min}) and (\ref{pr2}), we get (\ref{minr2}). 
From (\ref{dimSmin}) and  (\ref{minr2}), we obtain that 
 $$\dim (S) \leq  
 \left\lfloor \frac{r}{2} \right\rfloor
\left(n-  \left\lfloor \frac{r}{2} \right\rfloor\right)
  .$$ 
 
\end{proof}

Observe that in an analogous way 
we can prove  the following two theorems: 

\begin{thm} \label{signature} Let $p, q , n \in \mathbb{N}$ such that $ p+q \leq n$; then 
$$ a^{\mathbb{R}}_{sym}(n; p,q) \leq  \min{\{p,q\}} (n-  \min{\{p,q \}}
).$$
\end{thm}

\begin{proof}[Sketch of the proof]
Consider $S \in  {\cal A}_{sym}^{\mathbb{R}} (n; p \times q ) $. 
 We can suppose $S=D+ V$ where $D$ is the diagonal matrix whose diagonal is $(1, \ldots,1, -1, \ldots, -1, 0, \ldots, 0)$ where $1$ is repeated $p$ times and $-1$ is repeated $q$ times and $V$ 
  is a linear subspace of $S(n, \mathbb{R})$ and then argue as in the  proof of Theorem \ref{mio}.
\end{proof}

\begin{thm}\label{herm} (i) Let $n,r \in \mathbb{N}$ with $r \leq n$.
Then
$$a_{herm}(n;r) \leq 2 \left\lfloor \frac{r}{2} \right\rfloor (n-\left\lfloor \frac{r}{2} \right\rfloor) .$$ 

(ii)   Let $n,p, q \in \mathbb{N}$ with $p+ q  \leq n$.
Then
$$a_{herm}(n;p,q) \leq 2\min{\{p,q\}} (n-  \min{\{p,q \}}
).$$

\end{thm}

\begin{proof}[Sketch of the proof]
(i) 
 Let $R \in  {\cal A}_{herm} (n;r) $.  We can write $R$ as $ M +L$ where $M \in H(n)$ and $L$ is a linear subspace of $H(n)$. There exists a unitary matrix $U$ and a diagonal real matrix $P$ such that $ P \,
  {}^t \overline{U} M U P$ is 
   the 
 diagonal matrix $D$ whose diagonal is $(1, \ldots, 1 , -1, \ldots, -1, 0, \ldots, 0)$, where $1$ is repeated $p$ times for some $p$ and $-1$ is repeated $q$ times with $p+q =r$. Consider $S =  P\, {}^t \overline{U} R U P $; it  is 
 equal to $D+V$, where $V$ is a vector subspace of $H(n)$.

 Let $Z$ be the vector subspace of $H(n)$ generated by the matrices $E_{l,j}+E_{j,l}$ 
 and $i E_{l,j}-iE_{j,l}$ 
 for $l, j \in \{1, \ldots, p\}$ with $l < j$ and by the matrices $E_{l,l}$ for $l \in \{1, \ldots, p\}$.
 
 Let $U$ be the vector subspace of $H(n)$ generated by the matrices $E_{l,j}+E_{j,l}$ and $i E_{l,j}-iE_{j,l}$  for $l, j \in \{p+1, \ldots, r\}$ with $l < j$ and by the matrices $E_{l,l}$ for $l \in \{p+1, \ldots, r\}$.

 Let $W$ be the vector subspace of $H(n)$ generated by the matrices $E_{l,j}+E_{j,l}$ and $i E_{l,j}-iE_{j,l}$   for $l, j \in \{r+1, \ldots, n\}$ with $l < j$ and by the matrices $E_{l,l}$ for $l \in \{r+1, \ldots, n\}$.

Let $G$   be the vector subspace of $H(n)$
generated by the matrices $E_{l,j}+E_{j,l}$ 
and $iE_{l,j}-iE_{j,l}$ 
for $l \in \{r+1, \ldots, n\}$, $j\in \{p+1, \ldots , r\}$.

 As in the proof of Theorem \ref{mio}, we
 can  prove that 
$V \cap (Z+U+W+G)=\{0\}$. Hence 
$$ \dim (R)= \dim (S) = \dim ( V) \leq  \dim (H(n)) - \dim (Z+U+W+G)= 2p (n-p).$$
In an analogous way we can prove that $\dim (R) \leq 2 q (n-q)$. 

As in the proof of Theorem \ref{mio}
we can deduce that $\dim (R)
\leq   2 \left\lfloor \frac{r}{2} \right\rfloor (n-\left\lfloor \frac{r}{2} \right\rfloor) .$
 
The proof of (ii) is analogous.
\end{proof}

\begin{proof}[Proof of Theorem \ref{mio2}]
In order to prove that $a^{\mathbb{R}}(m \times n,r) $ is greater than or equal to 
$ r(n-r)+ \frac{r(r-1)}{2} $, i.e.  greater than or equal to 
$ rn- \frac{r(r+1)}{2} $,
consider the following affine subspace of $M(m \times n , \mathbb{R})$:
$$ S=\left\{ A \in M(m \times n , \mathbb{R}) \mid 
A_{i,i}=1 \; \forall i =1, \ldots, r , \; A_{i,j} =0 \; \forall (i,j) \; \mbox{\rm with } i>j  \; \mbox{\rm or }  i>r 
\right\}.$$
The dimension of $S$ is clearly $  r(n-r)+ \frac{r(r-1)}{2} $ and $S \in {\cal A}^{\mathbb{R}} ( m \times n;r)$ , so we get our inequality.

Now let us prove the other inequality.

Let $C \in  {\cal A}^{\mathbb{R}}(m \times n; r)$. 
 We want to prove that $\dim (C) \leq
 r(n-r)+ \frac{r(r-1)}{2} $. We can write $C$ as $A+ W$ where $A \in M(m \times n , \mathbb{R})$ and $W$ 
  is a linear subspace of $M(m \times n , \mathbb{R})$. Let $Q$ and $R$ be invertible matrices such that, if we denote  $Q^{-1} A R$ by $J$, we have that  
  $J_{i,i}=1$ for $i=1, \ldots, r$ and the other entries of $J$  are equal to zero.

 Let $V=  Q^{-1} W R $ and $S=   Q^{-1} C R =J+V$. Obviously $S  \in  {\cal A}^{\mathbb{R}}(m \times n; r)$; moreover, $ \dim(S)= \dim (C)$, so to prove that $\dim (C) \leq r(n-r)+ \frac{r(r-1)}{2}  $ it is sufficient to prove that $\dim (S) \leq r(n-r)+ \frac{r(r-1)}{2}  $.
 

Consider now the following subspaces of $M(m \times n, \mathbb{R})$: 
$$ Z=\{ A \in M(m \times n , \mathbb{R}) \mid A_{i,j}=0\; \forall (i,j)\; \mbox{\rm such that } i \neq j \;
\mbox{\rm and } (i \leq r \; \mbox{\rm or } j \leq r )\}$$ 
{\small
$$ T=\left\{ A \in M(m \times n , \mathbb{R}) \mid
\begin{array}{l}
A_{i,j}=0\; \forall (i,j)\; \mbox{\rm such that} \; i = j \;
\mbox{\rm or } ( i > r \; \mbox{\rm and } j > r)
 \;
\mbox{\rm or } j>m
; \\
A_{i,j}= A_{j,i} \; \forall (i,j) \; \mbox{\rm such that } j \leq m
\end{array}
\right\}.$$} 
We want to prove that 
$$ V \cap (Z+T)=\{0\}.$$
Let $\zeta  \in Z $ and $\tau \in T$ such that 
$\zeta+ \tau \in V$; we want to show that $ \zeta= \tau =0$.

We denote   $\zeta^{(r+1, \ldots , n)}_{(r+1, \ldots, m) }$ by $\zeta'$ and $\tau^{(1, \ldots ,r) }_{(r+1, \ldots, m)} $ by $\tau'$.

We consider four cases:

\underline{Case 1:  $\tau'=0$, $\zeta'=0 $.}

Since $\tau + \zeta \in V$  we have that $ J+ s(\tau + \zeta) $ must have rank $r$ for every $s \in \mathbb{R}$;  observe that  $ (\tau + \zeta) _{i,j}= 0 $ $\forall (i,j) $ with $i >r $ or $j >r$ 
(by the definition of $Z$ and $T$ and the fact that $\tau'=0$ and $ \zeta'=0$)
and that $ (\tau + \zeta)^{(1, \ldots, r)}_{(1, \ldots, r)}$ is symmetric; hence, by Lemma \ref{lemma2}, we can conclude that $ (\tau + \zeta)^{(1, \ldots, r)}_{(1, \ldots, r)}=0$ and then that   $\tau + \zeta=0$.

\underline{Case 2:  $\tau' = 0$, $\zeta' \neq 0$.}

Take $s\in \mathbb{R}-\{0\}$ such that $\det \left(I_r+ s  (\tau + \zeta)^{(1, \ldots, r)}_{(1, \ldots, r)} \right)$ is nonzero and $h \in \{r+1, \ldots, m\}$ and $l \in \{r+1, \ldots, n\}$  such that $\zeta_{h,l} \neq 0$.  Then, obviously, the matrix $(J + s  (\tau + \zeta) )^{(1, \ldots, r, l)}_{(1, \ldots, r, h)}$ is invertible, , which is impossible since $J + s  (\tau + \zeta) \in S \in {\cal A}^{\mathbb{R}}(m \times n; r)$.

\underline{Case 3:  there exists $h \in \{r+1, \ldots , m\}$ such that $\tau^{(h)} \neq 0$ and $\zeta_{(h,h)} = 0$.}

By Lemma \ref{lemma3},  we have that  $\det \left((J + s  (\tau + \zeta) )^{(1, \ldots, r, h)}_{(1, \ldots, r, h)}
\right)$ is 
a nonzero polynomial in $s$, so we can find $s$ 
such that $\det \left((J + s  (\tau + \zeta) )^{(1, \ldots, r, h)}_{(1, \ldots, r, h)} \right)\neq 0$, 
which is absurd since $J + s  (\tau + \zeta) \in S \in {\cal A}^{\mathbb{R}}(m \times n; r)$.

\underline{Case 4:  there exists $h \in \{r+1, \ldots , m\}$ such that $\tau^{(h)} \neq 0$ and $\zeta_{(h,h)} \neq 0$.}

 Observe that $\det \left( (J + s  (\tau + \zeta)^{(1, \ldots, r, h)}_{(1, \ldots, r, h)} \right)$ is a polynomial in $s$ with the term of  degree $0$ equal to $0$ and the coefficient  of the term of degree $1$ equal to $ \zeta_{h,h}$, which is nonzero; then there exists $s \in \mathbb{R}-\{0\}$ such that $\det \left((J + s  (\tau + \zeta) )^{(1, \ldots, r, h)}_{(1, \ldots, r, h)} \right)$ is nonzero; but this is impossible since $J + s  (\tau + \zeta) \in S \in {\cal A}^{\mathbb{R}}(m \times n; r)$.

Observe that the four cases we have considered are the only  possible ones because when $\tau'$ is nonzero we have one among Case 3 and Case 4. 
Thus we have proved that 
$ V \cap (Z+T)=\{0\}$.
Hence we have:
$$ \begin{array}{lll}
\dim (S) & = & \dim (V) \leq \dim M(m \times n, \mathbb{R})- \dim (Z+T) =  \vspace*{0.2cm} \\ 
& = & mn- \dim (Z)- \dim (T)=  r(n-r)+ \frac{r(r-1)}{2} 
\end{array}$$ 
and we can conclude.
\end{proof}

\begin{rem} Let $F[x_1,  \ldots, x_k ]$ denote the set of the polynomials in the indeterminates $x_1, \ldots, x_k$ with coefficients on a field $F$.
A matrix over $F[x_1,  \ldots, x_k ]$  is said an Affine Column Indipendent matrix, or ACI-matrix,
if its entries are polynomials of degree at most one and no indeterminate appears in two different columns. A completion of an ACI-matrix is an assignment of values in $F$ 
to the indeterminates $x_1, \ldots, x_k$;
 for instance, let us consider the matrix over $\mathbb{R}[x_1,\ldots, x_5]$ $$A =\begin{pmatrix} x_1 & x_3 & x_4+x_5 \\
 2x_1 +x_2 & -x_3 -1 & x_4-x_5 \\
 x_2+1 & 0 & 2x_4 \end{pmatrix}
;$$
it is an ACI matrix; if we assign the values $1,1,2,5,7$ respectively to $x_1, \ldots, x_5$, we get the completion of $A$
$$ \begin{pmatrix} 1 & 2 & 12 \\ 3 & -3 & -2 \\ 2 & 0 & 10 \end{pmatrix} .$$  

In \cite{H-Z} Huang and Zhan proved that all the completions of an $m \times n $ ACI-matrix $A$ over a field $F$ with $|F| \geq \max\{m,n+1\}$ have rank $r$ if and only if there exists a nonsingular constant $m \times m$ matrix $T$  and a permutation $n \times n $ matrix $Q$ 
such that $TA Q$ is equal to a matrix of the kind $$ \begin{pmatrix} B & \ast & \ast \\ 0 & 0 & \ast \\ 
0 & 0 & C \end{pmatrix} $$ for some ACI-matrices $B$ and $C$ which are square upper triangular with nonzero constant diagonal entries and whose orders sum to $r$.
Observe that the affine subspace given by  the matrices $ \begin{pmatrix} B & \ast & \ast \\ 0 & 0 & \ast \\ 
0 & 0 & C \end{pmatrix} $ with $B$ and $C$  square upper triangular ACI-matrices with nonzero constant diagonal entries and whose orders are respectively $k$ and  $r-k$ is equal to the following number:
$$ \frac{k(k-1)}{2}+ \frac{(r-k) (r-k-1)}{2} + k (n-k)+ (r-k)(m-k-r+k)= $$
$$= -\frac{r^2}{2} -\frac{r}{2} + k(n-m) +rm, $$
which obviously attains the maximum, i.e. $rn -\frac{r^2+ r}{2}$, when $k=r$. 
Let $M$ be a matrix over $F[x_1,  \ldots, x_k ]$ 
with the degree of every entry at most one.
Define $\tilde{M}$ to be the ACI-matrix 
obtained from $M$ in the following way: if an indeterminate $x_i$ appears in more than one column, say in the columns $j_1, \ldots, j_s$,    replace it with new indeterminates $x_i^{j_1},\ldots, x_i^{j_s}$, precisely replace $x_i$ in the $j_l$-th column with $x_i^{j_l}$ for $l=1, \ldots s$.

Observe that an affine subspace of $M(m \times n, F)$ corresponds to an $m \times n$ matrix over $F[x_1, \ldots, x_l]$  for some $l \in \mathbb{N}$ 
with the degree of every entry at most one and so an $m \times n$ ACI matrix such that all the completions have rank $r$ corresponds to an affine subspace of $M(m \times n, F)$ of constant rank $r$.

One might think that  it is  possible to deduce Theorem \ref{mio2}, in particular the inequality  
$a^{\mathbb{R}}(m \times n,r) \leq rn- \frac{r(r+1)}{2} $,
from Huang-Zhan`s result in the following way:
let $S \in {\cal A}^{\mathbb{R}}(m \times n; r)$ and consider 
the affine subspace $\tilde{S}$ given by the ``ACImade'' matrices of $S$, that is, given by 
 the matrices $\tilde{M}$ for $M \in S$; if $\tilde{S}$ were of constant rank, then it would correspond to an ACI matrix such that all the completions have rank $r$;
so   by Huang-Zhan`s result we would have $\dim (\tilde{S}) \leq  rn- \frac{r(r+1)}{2} $ and then $\dim (S) \leq  rn- \frac{r(r+1)}{2} $, but it is not true that 
the affine subspace $\tilde{S}$  is of constant rank
for any  affine subspace $S$ of  constant rank, 
 as the following example shows: let
 $$S =\left\{ \begin{pmatrix} 1 & s \\ s & -1  \end{pmatrix} \;| \; s \in \mathbb{R} 
 \right\};$$
 We have that $\tilde{S} =\left\{ \begin{pmatrix} 1 & s \\ t& -1  \end{pmatrix} \;| \; s ,t \in \mathbb{R} 
 \right\},$ which is not of constant rank.
\end{rem}

\begin{rem} \label{altricampi}
Observe that Theorem \ref{mio2} does not hold on every field $K$, as the following example shows:
consider the field $\mathbb{Z}/2$, $m=n=2$ and $r=1$; let 
$$ S= \begin{pmatrix}  1 & 0 \\ 0 & 0 \end{pmatrix} + \langle \begin{pmatrix} 1 & 0 \\ 0 & 1 \end{pmatrix} , \begin{pmatrix} 0 & 1 \\ 0 & 0 \end{pmatrix}  \rangle ;$$
the affine subspace $S$ is obviously of dimension $2$ and constant rank $1$, so  for $m=n=2$, $r=1$, we have that 
$a^{\mathbb{Z}/2}(m \times n,r) $ is different from $ rn- \frac{r(r+1)}{2} = 1$.
Anyway, the main problem to extend Theorem \ref{mio2} to other fields seems Lemma \ref{lemma2}, which we use in Case 1 of the proof. Precisely, observe that the argument to prove the inequality 
$a^{K}(m \times n,r)  \geq rn- \frac{r(r+1)}{2} $ works on any field $K$; as to the other inequality,
 it is easy to see that  the argument  in Cases 2,3,4 works for any field with cardinality greater than $r+2$ (in fact, such condition guarantees for any nonzero polynomial $p$ over $K$ in one variable of degree less than or equal to $r+1$ the existence of an element $s \in K-\{0\}$ such that $p(s) \neq 0$, in particular there exists $s\in K-\{0\}$ such that $\det \left(I_r+ s  (\tau + \zeta)^{(1, \ldots, r)}_{(1, \ldots, r)} \right) \neq 0 $  in Case 2 and there exists $s\in K-\{0\}$ such that  $\det \left((J + s  (\tau + \zeta) )^{(1, \ldots, r, h)}_{(1, \ldots, r, h)} \right)\neq 0$ in Cases 3 and 4); so the main problem to extend the theorem seems in Case 1, because we use Lemma \ref{lemma2}. 
 
Also for Theorem \ref{mio}, the main obstacle to extend the statement to other fields seems the necessity to extend Lemma \ref{lemma2}.

\end{rem}

{\bf Acknowledgments.}
This work was supported by the National Group for Algebraic and Geometric Structures, and their  Applications (GNSAGA-INdAM).

The author wishes to thank the anonymous referee for his/her comments, which helped to improve the paper.

\end{document}